\documentclass{elsart}               % The use of LaTeX2e is preferred.

\usepackage{graphicx}          % Include this line if your 
\usepackage{amssymb}    
\usepackage{amsmath}                           % you prefer.

\newtheorem{calc}[thm]{Calculation}
\newtheorem{example}[thm]{Example}
\renewcommand{\theenumi}{\alph{enumi}} \renewcommand{\theenumii}{\arabic{enumii}}

\begin{document}

\begin{frontmatter}

\title{On an implicit triangular decomposition of nonlinear control systems that are 1-flat - a constructive approach}

%\thanks[footnoteinfo]{This paper was not presented at any IFAC 
%meeting. Corresponding author M.~T.~Cicero. Tel. +XXXIX-VI-mmmxxi. 
%Fax +XXXIX-VI-mmmxxv.}

\author[JKU]{Markus Sch\"{o}berl}
\ead{markus.schoeberl@jku.at}
\author[JKU]{Kurt Schlacher}
\ead{kurt.schlacher@jku.at}
\address[JKU]{Johannes Kepler University Linz, Institute of Automatic Control and Control Systems Technology,
     Austria}
%\address[ACCM]{Austrian Center of Competence in Mechatronics (ACCM), Austria}
          
                             % keyword list or with the 
                                          % help of the Automatica 
                                          % keyword wizard

\begin{abstract}                          % Abstract of not more than 200 words.
We study the problem to provide a triangular form based on implicit differential equations for non-linear multi-input systems with respect to the flatness property. Furthermore, 
we suggest a constructive method for the transformation of a given system into that special triangular shape, if possible. 
The well known Brunovsky form, which is applicable with regard to the exact linearization problem, can be seen as special case of this implicit triangular form. 
A key tool in our investigation will be the construction of Cauchy characteristic vector fields that
additionally annihilate certain codistributions. In adapted coordinates this construction allows to single out variables whose time-evolution can be derived without any integration.   

\end{abstract}

\end{frontmatter}

\section{Introduction}

The concept of flatness introduced in \cite{Fliess,Fliessback} has
greatly influenced the control and systems theory community. The property
of a system to be flat allows for an elegant solution for many feed-forward
and/or feedback problems and is applicable for a big class of systems
including the linear and the nonlinear as well as the lumped- and
the distributed-parameter case. Within this paper we are interested
in the system class of nonlinear multi-input systems described by
ordinary differential equations. For this system class necessary and
sufficient conditions for flatness have been proposed in \cite{LevineFlatBook,LevineArt}
based on a polynomial matrix approach. Furthermore, for nonlinear
multi-input systems with special structure further results exist,
see e.g. \cite{Martin,MartinRouchon,PometESAIM,RathinamSIAM}.

Triangular forms are of special interest for nonlinear systems in
the context of exact linearization or flatness. Systems that are exactly
linearizable by static feedback can be converted to Brunovsky normal
form (a special case of the extended Goursat form, adapted to control
systems), see \cite{Tilbury}. It is well known that systems that
are flat but not exactly linearizable by static feedback can be transformed
into Brunovsky normal form only after a dynamic system extension (dynamic
compensator), see \cite{CharletSIAM}. For systems that are 0-flat
in \cite{Bououden} a nonlinear explicit triangular form has been
proposed. We consider systems that are 1-flat, the flat output may
depend on the state (0-flat) and the control (1-flat) but not on the
derivatives of the control. We will consider a\textit{ triangular
decomposition} based on \textit{implicit} differential equations,
and we will propose a constructive scheme how to transform a 1-flat
system into that special form, if possible (this gives rise to a sufficient
condition for a control system to be 1-flat). As in \cite{Tilbury},
where the extended Goursat form is discussed, we also make use of
the Pfaffian system representation such that 1-forms (covector-fields)
are used for the description of the (implicit) differential equations.
Furthermore, we will also make use of a filtration which is connected
to a triangular representation of implicit differential equations
in the case of flat systems. This is in contrast to the filtration
which is based on \textit{derived flags} used in the exact linearization
problem leading to a representation based on explicit differential
equations like the Brunovsky form, see again \cite{Tilbury}. It should
be noted that the \textit{implicit triangular decomposition }contains
the Brunovsky form as a special case. A different approach based on
differential forms, but associated with the tangent linear system,
can be found in \cite{Aranda}, where the so-called infinitesimal
Brunovsky form is considered.

This contribution can be seen as further developing the ideas presented
in \cite{SchlacherSchoeberlNolcos2007,SchoeberlNolcos2010}, where
a reduction and elimination procedure is considered to derive flat
outputs which (in contrast to this contribution) is not based on a
Pfaffian representation and the adequate tools from exterior algebra.
Preliminary results have been presented already in \cite{SchoeberlICCA},
where also an extended example, the VTOL, can be found, which has
been analyzed using different tools e.g. in \cite{Fliessback,LevineFlatBook}.

\textit{Notation:} Let $\mathcal{X}$ be an $n_{x}$ dimensional manifold,
equipped with local coordinates $x^{\alpha},\,\alpha=1,\ldots,n_{x}$,
i.e., $\mathrm{dim}(\mathcal{X})=\mathrm{dim}(x)=n_{x}$. We denote
the partial derivatives by $\partial_{x^{\alpha}}$ and $\partial_{x}=\partial_{x^{1}},\ldots,\partial_{x^{\alpha}},\ldots,\partial_{x^{n_{x}}}$.
We will make use of the Einstein convention on sums, namely $a_{i}b^{i}=\sum_{i=1}^{n}a_{i}b^{i}$
when the index range is clear from the context. Frequently, we will
use tensors in matrix representation together with index notation.
Given a matrix $m(x)$ the components are given as $m_{\alpha}^{j}(x)$
where the index $j$ corresponds to the rows and $\alpha$ to the
columns. We will need to numerate matrices such that we have for example
that $m_{3,\alpha}^{2,j}(x)$ are the components of a matrix $m_{3}^{2}(x)$.
If we multiply $m_{3}^{2}(x)$ with a vector $w=(w^{1},\ldots,w^{n_{w}})$
we have in components $m_{3}^{2}(x)w=m_{3,\alpha}^{2,j}(x)w^{\alpha}$
(where the components of $w$ are represented in a list, but in matrix
notation $w$ is interpreted as a column vector) and when a vector
$v$ is partitioned into blocks e.g. $v=(v^{1},v^{2})=(v^{1,1},\ldots,v^{1,n_{v^{1}}},v^{2,1},\ldots,v^{2,n_{v^{2}}})$
then we can compute e.g. $m_{3}^{2}(x)v^{1}=m_{3,\alpha}^{2,j}(x)v^{1,\alpha}$
(assuming appropriate dimensions of $m_{3}^{2}$ and $v^{1}$). We
will use the numeration of matrices and vectors to indicate to which
block they belong according to the implicit triangular form to be
defined.

\section{The triangular form}

Let us consider a nonlinear control system
\begin{equation}
\dot{x}=f(x,u)\label{eq:NLcontr}
\end{equation}
with $n_{x}$ states and $n_{u}$ independent inputs on a manifold
$\mathcal{X}\times\mathcal{U}$. Roughly speaking, the system (\ref{eq:NLcontr})
is flat ($(\kappa+1)-$flat), if there exist $n_{u}$ differentially
independent functions $y(x,u,\dot{u},\ldots,u^{(\kappa)})$, such
that the state $x$ and the control $u$ can be parameterized by $y$
and its successive time derivatives. Hence, flat systems enjoy the
characteristic feature that the (time) evolution of the state and
input (control) variables can be recovered from that of the flat output
without integration. The system is called 0-flat if $y$ depends solely
on $x$ and 1-flat if we have $y(x,u).$ For a rigorous definition
of differential flatness, see \cite{Fliess,LevineFlatBook}.

The Brunovsky-form, consisting of $n_{u}$ integrator chains, is the
most simple triangular structure that can be achieved for systems
(\ref{eq:NLcontr}) that are exactly linearizable by static feedback,
and hence 0-flat. A different (explicit) triangular form for 0-flat
systems has been proposed in \cite{Bououden} which is more general
than the Brunovsky form. To treat the case of 1-flat systems we will
present an \textit{implicit triangular form} for a special subclass
of systems of the form (\ref{eq:NLcontr}), that is useful regarding
the property to be 1-flat. 

The main idea in this contribution is to look for a diffeomorphism
$(x,u)=\varphi(z)$ such that in the new coordinates $z$ (corresponding
to all the system variables including the inputs) the system (\ref{eq:NLcontr})
can be represented in a special triangular shape consisting of implicit
differential equations, in a form such that the flat outputs can be
read off. We first will describe some properties of the coordinates
$z$, then we introduce the implicit triangular form (definition \ref{def:ImplTriag}),
and finally we discuss how to construct the map $\mbox{\ensuremath{\varphi}}$
in sections \ref{sec:Pfaffian-representation} and \ref{sec:A-constructive-algorithm}.

Therefore, let us consider a manifold $\mathcal{Z}$ where $\mathrm{dim}(\mathcal{Z})=n_{z}$
with coordinates $z$ which are partitioned in $m$ blocks of the
following form
\begin{equation}
\begin{array}{ccl}
z & = & (z^{1},z^{2},\ldots,z^{m})\\
 & = & ((z^{1,1},\ldots,z^{1,n_{z^{1}}}),\ldots,(z^{m,1},\ldots,z^{m,n_{z^{m}}}))
\end{array}\label{eq:z}
\end{equation}
where each $z^{i}$ consist of $n_{z^{i}}$ coordinates, i.e., $n_{z}=\sum_{i=1}^{m}n_{z^{i}}$.
In the following we will present a representation of the system (\ref{eq:NLcontr})
that among these $z$ coordinates we will find the flat outputs $y$
in the following form: 
\begin{equation}
(z^{1},z^{2},\ldots,z^{m})=(y^{1},(y^{2},\hat{z}^{2}),\ldots,(y^{m-1},\hat{z}^{m-1}),\hat{z}^{m})\label{eq:yzhat}
\end{equation}
with $n_{z^{i}}=n_{y^{i}}+n_{\hat{z}^{i}}$ and $n_{y^{j}}\geq0,\, j=2,\ldots,m-1$. 
\begin{rem}
The $n_{z^{i}}$ variables in the i-th block $z^{i}=(z^{i,1},\ldots,z^{i,n_{z^{i}}})$
are decomposed into $n_{z^{i}}=n_{y^{i}}+n_{\hat{z}^{i}}$ variables
according to $z^{i}=(y^{i},\hat{z}^{i})$. In the forthcoming we will
show that $y^{i}$ (possibly empty for $i>1$) will be part of the
flat output and the $\hat{z}^{i}$ will be called non-derivative variables
(since they appear non-differentiated in certain blocks of the triangular
form). 
\end{rem}
A desirable structure to study 1-flat system is defined by the following
implicit differential equations which are decomposed into $n_{b}$
blocks.
\begin{defn}
\label{def:ImplTriag}The implicit differential equations $\Xi_{e}^{i}=0\,,\, i=1,\ldots,n_{b}$
and $m=n_{b}+1$ based on the partition (\ref{eq:z}) and (\ref{eq:yzhat})
given as 
\begin{eqnarray}
\Xi_{e}^{1} & : & a_{1,\alpha_{1}}^{1,j_{1}}\dot{z}^{1,\alpha_{1}}-b^{1,j_{1}}\nonumber \\
\Xi_{e}^{2} & : & a_{1,\alpha_{1}}^{2,j_{2}}\dot{z}^{1,\alpha_{1}}+a_{2,\alpha_{2}}^{2,j_{2}}\dot{z}^{2,\alpha_{2}}-b^{2,j_{2}}\nonumber \\
 &  & \vdots\label{eq:Triag1}\\
\Xi_{e}^{n_{b}} & : & a_{1,\alpha_{1}}^{n_{b},j_{n_{b}}}\dot{z}^{1,\alpha_{1}}+\ldots+a_{n_{b},\alpha_{n_{b}}}^{n_{b},j_{n_{b}}}\dot{z}^{n_{b},\alpha_{n_{b}}}-b^{n_{b},j_{n_{b}}}\nonumber 
\end{eqnarray}
are termed the implicit triangular form with $j_{i}=1,\ldots,\mathrm{dim}(\Xi_{e}^{i})$
and $\alpha_{i}=1,\ldots,n_{z^{i}}$ for $i=1,\ldots,n_{b}$, which
possesses the following properties
\begin{enumerate}
\item the matrices $a_{k}^{i}$ and (the vectors) $b^{i}$ meet 
\begin{eqnarray}
a_{k,\alpha_{k}}^{i,j_{i}} & = & a_{k,\alpha_{k}}^{i,j_{i}}(z^{1},\ldots,z^{i},\hat{z}^{i+1})\,,\,\label{eq:ab}\\
b^{i,j_{i}} & = & b^{i,j_{i}}(z^{1},\ldots,z^{i},\hat{z}^{i+1})\nonumber 
\end{eqnarray}

\item $\mathrm{dim}(\Xi_{e}^{i})=\mathrm{dim}(\hat{z}^{i+1})$, and the
Jacobian matrices $[\partial_{\hat{z}^{i+1}}\Xi_{e}^{i}]$ are regular
for all $i=1,\ldots,n_{b}=m-1$.
\end{enumerate}
\end{defn}

It should be noted that $j_{i},\alpha_{k}$ in (\ref{eq:Triag1})
are indices corresponding to the rows and the columns of the matrices
$a_{k}^{i}$ and vectors $b^{i}$ respectively (summation over the
$\alpha_{i}$), where each subsystem can be represented as $\Xi_{e}^{i}:a_{k,\alpha_{k}}^{i,j_{i}}\dot{z}^{k,\alpha_{k}}-b^{i,j_{i}}$
(summation over $k$ and $\alpha_{k}$ with $k\leq i$) and that due
to (\ref{eq:ab}) the dependence on the $z$ coordinates is arranged
in a triangular manner, as demonstrated in the following example.
\begin{example}
\label{ex2}A system in triangular form with $n_{b}=3$ (3 blocks
for the equations) and thus $m=n_{b}+1=4$ (4 blocks in $z$) in matrix
notation reads as
\begin{eqnarray*}
\Xi_{e}^{1} & : & a_{1}^{1}(z^{1},\hat{z}^{2})\dot{z}^{1}-b^{1}(z^{1},\hat{z}^{2})\\
\Xi_{e}^{2} & : & a_{1}^{2}(z^{1},z^{2},\hat{z}^{3})\dot{z}^{1}+a_{2}^{2}(z^{1},z^{2},\hat{z}^{3})\dot{z}^{2}-b^{2}(z^{1},z^{2},\hat{z}^{3})\\
\Xi_{e}^{3} & : & {\displaystyle a_{1}^{3}}(z^{1},z^{2},z^{3},\hat{z}^{4})\dot{z}^{1}+\ldots+{\displaystyle a_{3}^{3}}(z^{1},z^{2},z^{3},\hat{z}^{4})\dot{z}^{3}-b^{3}
\end{eqnarray*}
where $b^{3}=b^{3}(z^{1},z^{2},z^{3},\hat{z}^{4})$ and 
\[
z=(z^{1},z^{2},z^{3},z^{4})=(y^{1},(y^{2},\hat{z}^{2}),(y^{3},\hat{z}^{3}),\hat{z}^{4})
\]
such that $y^{2}$ and/or $y^{3}$ are possibly empty (but they need
not, $n_{y^{2}},n_{y^{3}}\geq0$ ). We require $\mathrm{dim}(\Xi_{e}^{i})=\mathrm{dim}(\hat{z}^{i+1})$,
and the Jacobian matrices $[\partial_{\hat{z}^{i+1}}\Xi_{e}^{i}]$
are regular for all $i=1,\ldots,3$ such that $\hat{z}^{i+1}$ can
be computed by means of the implicit function theorem from $\Xi_{e}^{i}.$
\end{example}

\begin{lem}
\label{lem:FlatLemma} $y$ is a flat output for the system (\ref{eq:Triag1}). 
\end{lem}
To prove this Lemma, we consider the implicit equations $\Xi_{e}^{1}=0$.
Then $z^{1}(t)=y^{1}(t)$ can be assigned freely, and $\hat{z}^{2}(t)$
can be computed, where we make use of the implicit function theorem.
We continue with the equations $\Xi_{e}^{2}=0$. Two scenarios are
possible: $z^{2}=\hat{z}^{2}$, then it can be easily checked that
$\hat{z}^{3}(t)$ can be computed using the same argument as for $\hat{z}^{2}(t)$,
since $z^{1}(t)$ and $z^{2}(t)$ are already given. If $z^{2}=(y^{2},\hat{z}^{2})$
then $y^{2}(t)$ can be chosen freely, since the rank criteria is
met for $\hat{z}^{3}(t),$ which again can be computed. By continuing
this procedure, we end up with the equations $\Xi_{e}^{n_{b}}=0$
from which $\hat{z}^{m}(t)$ can be computed since at this stage $z^{1}(t),\ldots,z^{m-1}(t)$
are already known. This clearly shows that $(y^{1},\ldots,y^{m-1})$
is a flat output for (\ref{eq:Triag1}).
\begin{prop}
(a sufficient condition)\label{cor:The-control-system} The control
system (\ref{eq:NLcontr}) is 1-flat if we can find locally a diffeomorphism
$(x,u)=\varphi(z^{1},\ldots,z^{m})$ with $\sum_{i=1}^{m}n_{z^{i}}=n_{u}+n_{x}$,
such that it can be represented in the form (\ref{eq:Triag1}) and
$\sum_{i=1}^{m-1}n_{y^{i}}=n_{u}$.
\end{prop}
The fact that this proposition is sufficient for 1-flat systems, comes
from the observation that the flat outputs are among the $z$ coordinates,
and therefore clearly a function of $x$ and $u$ - as a special case
also 0-flat systems are included. Furthermore, it should be noted
that we consider a diffeomorphism which implies that we do not increase
the dimension of the system variables.

The goal is now to provide a constructive algorithm that transforms
a nonlinear multi-input system (\ref{eq:NLcontr}), if possible, into
the form (\ref{eq:Triag1}). Before we will analyze this in detail,
let us consider an example.

\subsection{A motivating example}

We consider a system with three state variables $(x^{1},x^{2},x^{3})$
and two control inputs $(u^{1},u^{2})$ of the form 

\begin{equation}
\dot{x}^{1}=u^{1},\,\,\,\dot{x}^{2}=u^{2}\,,\,\,\dot{x}^{3}=\sin\left(\frac{u^{1}}{u^{2}}\right)\label{eq:ExSin}
\end{equation}
also analyzed in \cite{LevineFlatBook,SchlacherSchoeberlNolcos2007}
using a different approach. Let us introduce the local coordinate
transformation $(x^{1},x^{2},x^{3},u^{1},u^{2})=\varphi(y^{1},\hat{z}^{2},y^{2},\hat{z}^{3},\hat{z}^{4})$
together with its inverse 
\begin{equation}
\begin{array}{cllccccccccccl}
x^{1} & = & \hat{z}^{2}\hat{z}^{3} &  &  &  &  &  &  &  &  & y^{1} & = & x^{3}\\
x^{2} & = & \hat{z}^{3}+y^{2} &  &  &  &  &  &  &  &  & \hat{z}^{2} & = & \frac{u^{1}}{u^{2}}\\
x^{3} & = & y^{1} &  &  &  &  &  &  &  &  & y^{2} & = & x^{2}-x^{1}\frac{u^{2}}{u^{1}}\\
u^{1} & = & e^{\hat{z}^{4}}\hat{z}^{2} &  &  &  &  &  &  &  &  & \hat{z}^{3} & = & x^{1}\frac{u^{2}}{u^{1}}\\
u^{2} & = & e^{\hat{z}^{4}} &  &  &  &  &  &  &  &  & \hat{z}^{4} & = & \mathrm{ln}(u^{2}).
\end{array}\label{eq:Trafo1}
\end{equation}
Then the system (\ref{eq:ExSin}) in the new coordinates can be represented
as
\begin{equation}
\begin{array}{ccc}
\dot{y}^{1}-\sin\left(\hat{z}^{2}\right) & = & 0\\
-\dot{y}^{2}\hat{z}^{2}+\dot{\hat{z}}^{2}\hat{z}^{3} & = & 0\\
\dot{y}^{2}+\dot{\hat{z}}^{3}-e^{\hat{z}^{4}} & = & 0
\end{array}\label{eq:SinNormal}
\end{equation}
(by a suitable combination of the equations) which is an implicit
system of differential equations. Following the proof of Lemma \ref{lem:FlatLemma}
it can be seen that the flat outputs are obviously $y^{1}$ and $y^{2}$
and in $(x,u)$ coordinates they read as $y^{1}=x^{3},\, y^{2}=x^{2}-x^{1}\frac{u^{2}}{u^{1}}$
based on (\ref{eq:Trafo1}). 
\begin{example}
The system (\ref{eq:SinNormal}) possesses the structure (\ref{eq:Triag1})
as in example \ref{ex2} with $n_{b}=3,\, m=4$ and $\mathrm{dim}(z^{1})=\mathrm{dim}(z^{3})=\mathrm{dim}(z^{4})=1$,
$\mathrm{dim}(z^{2})=2$ and the matrices
\[
\begin{array}{cccccccccccccccl}
a_{1}^{1} & = & 1 & , &  &  &  &  &  &  &  &  &  & b^{1} & = & \sin\left(\hat{z}^{2}\right)\\
a_{1}^{2} & = & 0 & , & a_{2}^{2} & = & [-\hat{z}^{2},\hat{z}^{3}] & , &  &  &  &  &  & b^{2} & = & 0\\
a_{1}^{3} & = & 0 & , & a_{2}^{3} & = & [1\,\,\,,\,\,0\,] & , & a_{3}^{3} & =1 & , &  &  & b^{3} & = & e^{\hat{z}^{4}}.
\end{array}
\]

\end{example}
The key question is now, how to derive the coordinate transformation
(\ref{eq:Trafo1}) and how must the equations be combined, such that
the form (\ref{eq:SinNormal}) can be obtained. These questions will
both be answered at once by using a Pfaffian system representation.

\section{Pfaffian representation\label{sec:Pfaffian-representation}}

We will use tools from exterior algebra and Pfaffian systems in the
sequel where we refer for detailed information to \cite{Bryant} and
references therein. For a representation of nonlinear control systems
in a Pfaffian form (with the focus on exact linearization with static
feedback), see e.g. \cite{Tilbury} and references therein. It should
be noted that we do not base our considerations on the tangent linear
system, as it is used for instance in \cite{Aranda}.

\subsection{Exterior Algebra and Properties of Pfaffian systems}

We denote by $\mathrm{d}\omega$ the exterior derivative of the $k$-form
$\omega$ and by $v\rfloor\omega$ the contraction (interior product)
of $\omega$ by the vector field $v$. The exterior product (wedge
product) is denoted by $\wedge$. 

A Pfaffian system $P$ on an $n_{\zeta}$-dimensional manifold $\mathfrak{Z}$
with coordinates $(\zeta^{\alpha})$ $\alpha=1,\ldots,n_{\zeta}$
can be identified with a codistribution $P=\{\omega^{1},\ldots,\omega^{n_{P}}\}$,
with $\omega^{i}=m_{\alpha}^{i}(\zeta)\mathrm{d}\zeta^{\alpha}$.
The\textit{ annihilator} of a Pfaffian system $P$ is a distribution
on $\mathfrak{Z}$ denoted by 
\[
P^{\perp}:=\left\{ w\in\mathcal{T}(\mathfrak{Z}),\, w\rfloor\omega=0,\,\forall\omega\in P\right\} .
\]
The \textit{derived flag }of the Pfaffian System $P$ is the descending
chain of Pfaffian systems $P^{(0)}\supset P^{(1)}\supset P^{(2)}\supset\ldots$
with $P^{(0)}=P$ and
\[
P^{(k+1)}:=\{\omega\in P^{(k)},\mathrm{d}\omega=0\,\mathrm{mod}\, P^{(k)}\}.
\]
\textit{Cauchy characteristic vector fields} $v$ of $P$ meet 
\begin{equation}
v\rfloor P=0\,,\,\,\,\,\, v\rfloor\mathrm{d}P\subset P.\label{eq:Cauchy1}
\end{equation}
The importance of Cauchy characteristic vector fields lies in the
fact that the Pfaffian system $P$ can be written using $n_{\zeta}-c$
coordinates (after a suitable coordinate transformation), where $c$
denotes the number of all independent Cauchy characteristic vector
fields. The distribution formed by the Cauchy characteristic vector
fields is denoted by $\mathcal{C}(P)$ which is involutive by construction,
see \cite{Bryant}. The desired coordinates can be constructed by
means of the Straightening out theorem (Frobenius theorem) such that
coordinates are introduced that are adapted to the involutive distribution
$\mathcal{C}(P)$. Indeed, by choosing adapted coordinates $(\bar{\zeta},\tilde{\zeta})=(\bar{\zeta}^{1},\ldots,\bar{\zeta}^{c},\tilde{\zeta}^{c+1},\ldots,\tilde{\zeta}^{n_{\zeta}})$
such that $ $$\mathcal{C}(P)=\{\partial_{\bar{\zeta}}\}$ the Pfaffian
system can be represented solely using the $\tilde{\zeta}$ coordinates,
see again \cite{Bryant}. 

In the following we consider \textit{time-invariant} dynamical systems
represented as\textit{ Pfaffian systems} on bundles. These systems
possess a fibration with respect to the time-manifold, i.e. let $\mathfrak{X}$
be an $n_{\xi}$-dimensional manifold, the corresponding fibration
is $\mathfrak{X}\times\mathbb{R}\rightarrow\mathbb{R}$, then a time-invariant
Pfaffian system $P$ is identified with a codistribution on the $(n_{\xi}+1)$-dimensional
manifold $\mathfrak{X}\times\mathbb{R}$ and is locally spanned by
1-forms $\omega$ of the form
\begin{equation}
\omega^{i}=m_{\alpha}^{i}(\xi)\mathrm{d}\xi^{\alpha}-n^{i}(\xi)\mathrm{d}t\,,\,\,\, P=\{\omega^{1},\ldots,\omega^{n_{P}}\}\label{eq:PfaffStandard}
\end{equation}
with $n_{P}=\mathrm{dim}(P)$. To $\omega^{i}$ there correspond the
implicit differential equations $\omega_{e}^{i}=0$ with $\omega_{e}^{i}=m_{\alpha}^{i}(\xi)\dot{\xi}^{\alpha}-n^{i}(\xi)$,
where $\dot{\xi}$ denotes the time-derivative.

Due to the fibration with respect to the time manifold we can introduce
a special kind of annihilator.
\begin{defn}
The vertical annihilator of (\ref{eq:PfaffStandard}) denoted by $\mathcal{V}(P)^{\perp}$
is defined to be the annihilator of the extended Pfaffian system $\{P,\mathrm{d}t\}$. 
\end{defn}
It is clear that $\mathcal{V}(P)^{\perp}\subset P{}^{\perp}$, i.e.,
one picks only those vector fields in $P{}^{\perp}$ which are tangential
to the fibration (those that do not include a $\partial_{t}$ component).
\begin{example}
Let us consider the explicit control system $\dot{x}=f(x,u)$ written
as a Pfaffian system $P=\{\omega^{\alpha_{x}}\}$ on a manifold with
coordinates $(t,x,u)$ which is fibred over the time with 
\[
\omega^{\alpha_{x}}=\mathrm{d}x^{\alpha_{x}}-f^{\alpha_{x}}(x,u)\mathrm{d}t\,,\,\,\,\,\alpha_{x}=1,\ldots,n_{x}
\]
i.e., $\xi=(x,u)$. Then we have $P^{\perp}=\{\partial_{t}+f^{\alpha_{x}}(x,u)\partial_{x^{\alpha_{x}}},\partial_{u}\}$
as well as $\mathcal{V}(P)^{\perp}=\{\partial_{u}\}$.
\end{example}

\begin{defn}
We call a Pfaffian system $P$ as in (\ref{eq:PfaffStandard}) \textit{parameterizable}
with respect to $\hat{\xi}$ when we can find appropriate coordinates
$(\bar{\xi},\hat{\xi})$ as well as a diffeomorphism $\xi=\psi(\bar{\xi},\hat{\xi})$
with $n_{\xi}=n_{\bar{\xi}}+n_{\hat{\xi}}$ where $n_{\hat{\xi}}=n_{P}=\mathrm{dim}(P)$
such that $P$ is represented as 
\begin{equation}
\bar{\omega}^{i}=\psi^{*}(\omega^{i})=q_{\alpha}^{i}(\bar{\xi},\hat{\xi})\mathrm{d}\bar{\xi}{}^{\alpha}-r^{i}(\bar{\xi},\hat{\xi})\mathrm{d}t\label{eq:PfaffNonDer}
\end{equation}
$\alpha=1,\ldots,n_{\bar{\xi}}$ and such that the differential equations
corresponding to (\ref{eq:PfaffNonDer}), i.e. $\bar{\omega}_{e}^{i}=0$
with $\bar{\omega}_{e}^{i}=q_{\alpha}^{i}(\bar{\xi},\hat{\xi})\dot{\bar{\xi}}^{\alpha}-r^{i}(\bar{\xi},\hat{\xi})$,
fulfill that the Jacobian matrix $[\partial_{\hat{\xi}}\bar{\omega}_{e}^{i}]$
is regular (and quadratic since $n_{\hat{\xi}}=n_{P}$). 
\end{defn}
The variables $\hat{\xi}$ are termed \textit{non-derivative} variables
and by the implicit function theorem we locally have $\hat{\xi}=g(\bar{\xi},\dot{\bar{\xi}})$.
Furthermore, it should be noted that $\{\partial_{\hat{\xi}}\}\subset\mathcal{V}(P)^{\bot}$
holds (since no $\mathrm{d}\hat{\xi}$ appears).

\subsection{The implicit triangular form in Pfaffian representation }

Let us consider the implicit differential equations as in (\ref{eq:Triag1})
written using differential forms on the bundle $\mathcal{Z}\times\mathbb{R}\rightarrow\mathbb{R}$
with the same properties as described in definition \ref{def:ImplTriag}
where $\Xi^{i}$ corresponds to the Pfaffian representation of $\Xi_{e}^{i}$
\begin{eqnarray}
\Xi^{1} & : & a_{1,\alpha_{1}}^{1,j_{1}}\mathrm{d}z^{1,\alpha_{1}}-b^{1,j_{1}}\mathrm{d}t\nonumber \\
\Xi^{2} & : & a_{1,\alpha_{1}}^{2,j_{2}}\mathrm{d}z^{1,\alpha_{1}}+a_{2,\alpha_{2}}^{2,j_{2}}\mathrm{d}z^{2,\alpha_{2}}-b^{2,j_{2}}\mathrm{d}t\nonumber \\
 &  & \vdots\label{eq:TRiagPfaff}\\
\Xi^{n_{b}} & : & {\displaystyle a_{1,\alpha_{1}}^{n_{b},j_{n_{b}}}}\mathrm{d}z^{1,\alpha_{1}}+\ldots+{\displaystyle a_{n_{b},\alpha_{n_{b}}}^{n_{b},j_{n_{b}}}}\mathrm{d}z^{n_{b},\alpha_{n_{b}}}-b^{n_{b},j_{n_{b}}}\mathrm{d}t\nonumber 
\end{eqnarray}
Let us denote by%
\footnote{The subscript $d$ will always refer to a representation based on
the desired triangular decomposition (\ref{eq:TRiagPfaff}).%
} $S_{d,0}$ the system (\ref{eq:TRiagPfaff}) and by $S_{d,k}=\{\Xi^{1},\ldots,\Xi^{n_{b}-k}\}.$ 
\begin{prop}
\label{prop:Properties}The system (\ref{eq:TRiagPfaff}) with $m=n_{b}+1$
enjoys the following properties 
\begin{enumerate}
\item $\{\partial_{\hat{z}^{m-k}}\}\subset\mathcal{V}(S_{d,k})^{\perp}$
are involutive distributions and $\partial_{\hat{z}^{m-k}}\in\mathcal{C}(S_{d,k+1})$
for $k=0,\ldots,n_{b}-1.$
\item Each subsystem $\Xi^{k}$ is parameterizable with respect to the non-derivative
variable $\hat{z}^{k+1}$, i.e. 
\[
\hat{z}^{k+1}=g^{k+1}(z^{1},\ldots,z^{k},\dot{z}^{1},\ldots,\dot{z}^{k})
\]
for $k=1,\ldots,n_{b}.$
\item If in $\Xi^{k}$, $z^{k}=(y^{k},\hat{z}^{k})$ such that $n_{z^{k}}>n_{\hat{z}^{k}}$,
i.e., variables $y^{k}$ are present, then $\partial_{y^{k}}\in\mathcal{C}(S_{d,m-k}).$
\end{enumerate}
\end{prop}
The proof of this proposition is straightforward and follows from
the structure of (\ref{eq:TRiagPfaff}) together with the special
structure of the $a_{k,\alpha_{k}}^{i,j_{i}}$ and $b^{i,j_{i}}$
according to (\ref{eq:ab}) as in definition \ref{def:ImplTriag}.
\begin{cor}
The implicit triangular decomposition (\ref{eq:TRiagPfaff}) gives
rise to the decomposition of $S_{d,0}$ into a sequence of Pfaffian
systems 
\begin{equation}
\ldots\subset S_{d,2}\subset S_{d,1}\subset S_{d,0}\label{eq:Seq1}
\end{equation}
as well as to splittings of the form $S_{d,i}=S_{d,i+1}\oplus S_{d,i+1,c}$,
where all the $S_{d,i+1,c}$ are parameterizable with respect to the
corresponding non-derivative variables $\hat{z}$.
\end{cor}

\begin{example}
(Example \ref{ex2} cont.) Following the notations in proposition
\ref{prop:Properties} we have $S_{d,0}=\{\Xi^{1},\Xi^{2},\Xi^{3}\}$,
$S_{d,1}=\{\Xi^{1},\Xi^{2}\}$ and $S_{d,2}=\{\Xi^{1}\}$ since $n_{b}=3$.
We observe that $\{\partial_{\hat{z}^{4}}\}\subset\mathcal{V}(S_{d,0})^{\perp}$
as well as $\{\partial_{\hat{z}^{4}}\}\subset\mathcal{C}(S_{d,1})$
which is obvious since $\hat{z}^{4}$ are non-derivative variables
which only appear in $\Xi^{3}$. The same holds true regarding $\hat{z}^{3}$
where now $\{\partial_{\hat{z}^{3}}\}\subset\mathcal{V}(S_{d,1})^{\perp}$
and $\{\partial_{\hat{z}^{3}}\}\subset\mathcal{C}(S_{d,2})$. Proposition
\ref{prop:Properties} (c) means for instance that if $z^{2}=(y^{2},\hat{z}^{2})$
then $\partial_{y^{2}}\subset\mathcal{C}(S_{d,2})$ since in $\Xi^{1}$
only $\hat{z}^{2}$ appears.
\end{example}

\section{A constructive algorithm\label{sec:A-constructive-algorithm}}

The goal is now to develop a constructive scheme that subsequently
creates this sequence (\ref{eq:Seq1}) as well as appropriate coordinate
transformations based on a given control system of the form $S_{0}=\{\omega_{0}^{\alpha_{x}}\}$
\begin{equation}
\omega_{0}^{\alpha_{x}}=\mathrm{d}x^{\alpha_{x}}-f^{\alpha_{x}}(x,u)\mathrm{d}t.\label{eq:Pfaff1Org}
\end{equation}
The starting point of the scheme is the explicit system $S_{0}$ but
since linear combinations of the $\omega_{0}^{\alpha_{x}}$ lead to
implicit equations in general we demonstrate the constructive method
with the system $S_{k}=\{\omega_{k}^{i}\}$ (here the index $k$ refers
to the $k-th$ iteration of the reduction process) with 
\begin{equation}
\omega_{k}^{i}=m_{\alpha}^{i}(\xi)\mathrm{d}\xi^{\alpha}-n^{i}(\xi)\mathrm{d}t\label{eq:Pfaff1OrgImp}
\end{equation}
with $i=1,\ldots,n_{e}$ and $n_{\xi}>n_{e}$, where we denote by
$\xi$ all the system variables. (\ref{eq:Pfaff1Org}) is a special
case of (\ref{eq:Pfaff1OrgImp}), i.e. $\xi=(x,u)$ in $S_{0}$. The
following steps need to be performed 
\begin{enumerate}
\item Computation of $\mathcal{V}(S_{k})^{\perp}$, since these elements
correspond to non-derivative variables. Choosing of an involutive
$\mathcal{F}_{k}\subset\mathcal{V}(S_{k})^{\perp}$ corresponds to
a selection of non-derivative variables called $\hat{w}_{k}$. (This
correspondence becomes obvious in an adapted coordinate chart to be
constructed by means of the Straightening out theorem.)
\item Construction of a splitting $S_{k}=S_{k+1}\oplus S_{k+1,c}$ such
that $\mathcal{F}_{k}\subset\mathcal{C}(S_{k+1})$, since this guarantees
that $S_{k+1}$ is independent of $\hat{w}_{k}$. 
\item Check, if $S_{k+1,c}$ is parameterizable with respect to the $\hat{w}_{k}$,
which is possible only if $\mathrm{dim}(S_{k})=\mathrm{dim}(S_{k+1})+\mathrm{dim}(\mathcal{F}_{k})$
holds.
\end{enumerate}
The whole procedure will then be continued with $S_{k+1}.$

\subsection{The k-th step of the system decomposition}

The constructive scheme rests on the following proposition
\begin{prop}
\label{prop:w_what}Let us consider the system $S_{k}=\{\omega_{k}^{i}\}$
with $\omega_{k}^{i}$ as in (\ref{eq:Pfaff1OrgImp}). If we find
an involutive distribution $\mathcal{F}_{k}$ with $\mathcal{F}_{k}\subset\mathcal{V}(S_{k})^{\perp}$
and a sub-codistribution $S_{k+1}\subset S_{k}$ such that $\mathcal{F}_{k}\subset\mathcal{C}(S_{k+1})$
is met, then we obtain a splitting $S_{k}=S_{k+1}\oplus S_{k+1,c}$
with 
\begin{equation}
\begin{array}{ccl}
\text{\ensuremath{S}}_{k+1} & : & \omega_{k+1}^{i}=a_{\alpha}^{i}(w_{k})\mathrm{d}w_{k}{}^{\alpha}-b^{i}(w_{k})\mathrm{d}t\\
S_{k+1,c} & : & \omega_{k+1,c}^{j}=a_{\alpha,c}^{j}(w_{k},\hat{w}_{k})\mathrm{d}w_{k}{}^{\alpha}-b_{c}^{j}(w_{k},\hat{w}_{k})\mathrm{d}t
\end{array}\label{eq:Decomp1}
\end{equation}
in adapted coordinates $(w_{k},\hat{w_{k}})$ by using a diffeomorphism
$\xi=\varphi_{k}(w_{k},\hat{w}_{k})$ with $n_{\xi}=n_{w_{k}}+n_{\hat{w}_{k}}.$ 
\end{prop}
The adapted coordinates can be constructed by means of the Straightening
out theorem since $\mathcal{F}_{k}$ is involutive, such that in new
coordinates $\mathcal{F}_{k}=\{\partial_{\hat{w}_{k}}\}$. In these
adapted coordinates $\partial_{\hat{w}_{k}}\subset\mathcal{V}(S_{k})^{\bot}$
as well as $\partial_{\hat{w}_{k}}\subset\mathcal{C}(S_{k+1})$ is
met, therefore no $\mathrm{d}\hat{w}_{k}$ can appear and a basis
of $S_{k+1}$ must exist which is independent of the $\hat{w}_{k}$
coordinates, since $\partial_{\hat{w}_{k}}\subset\mathcal{C}(S_{k+1})$.
Furthermore, if the system $S_{k+1,c}$ is parameterizable with respect
to $\hat{w}_{k}$ and if the system $S_{k+1}$ that can be expressed
in the $w_{k}$ coordinates possesses a non-trivial Cauchy characteristic,
then these redundant variables are candidates for possible flat outputs.
Based on these considerations we state the following corollary which
additionally includes the parameterization criteria, such that proposition
(\ref{prop:w_what}) is connected with the triangular form (\ref{eq:Triag1}),
respectively (\ref{eq:TRiagPfaff}).
\begin{cor}
\label{cor:SequenceOrig}The system $S_{0}$ (\ref{eq:Pfaff1Org})
can be transformed into the form (\ref{eq:TRiagPfaff}) if we find
a sequence of codistributions
\[
\ldots\subset S_{2}\subset S_{1}\subset S_{0}
\]
as well as involutive distributions $\mathcal{F}_{l}$ that meet $\mathcal{F}_{l}\subset\mathcal{V}(S_{l})^{\perp}$
as well as $\mathcal{F}_{l}\subset\mathcal{C}(S_{l+1})$ for $l\geq0$
such that the systems $S_{l+1,c}$ according to $S_{l}=S_{l+1}\oplus S_{l+1,c}$
are parameterizable with respect to $\mathcal{F}_{l}$. Then also
$\mathrm{dim}(S_{l+1,c})=\mathrm{dim}(\mathcal{F}_{l})$ holds where
we assume that each $S_{i}$ is represented by a minimal number of
variables. 
\end{cor}
This sequence ends when we have a decomposition of the form $S_{k^{*}}=S_{k^{*}+1}\oplus S_{k^{*}+1,c}$
with $S_{k^{*}+1}$ the empty system, which means that $S_{k^{*}}$
is a parameterizable system. This iterative scheme has therefore to
be continued until a parameterizable system is obtained. It should
be noted that in practice the effective computation of $\mathcal{F}_{l}$
and $S_{l+1}$ such that additionally parametrization is guaranteed
for all elements of the sequence is a difficult task. We will comment
on computational issues in section \ref{sub:A-constructive-method}
and demonstrate on an example a possible strategy.

\subsection{The connection with the derived flag}

In this short paragraph we want to discuss how the derived flag, see
\cite{Bryant} and its application to the exact linearizability problem
as described e.g. in \cite{Tilbury} is connected to our filtration,
as in corollary \ref{cor:SequenceOrig}.

Let us introduce an adapted basis for $S_{k},\,\mathrm{dim}(S_{k})=h$
(with respect to the derived flag) which is $S_{k}=\{\bar{\Omega}_{k}^{1},\ldots,\bar{\Omega}_{k}^{j},\Omega_{k}^{j+1},\ldots,\Omega_{k}^{h}\}$.
The first derived flag of $S_{k}$, denoted by $S_{k}^{(1)}$, meets
$\bar{\Omega}_{k}^{j}\in S_{k}^{(1)}\subset S_{k}\,,\, j=1,\ldots,\mathrm{dim}(S_{k}^{(1)})$
such that
\begin{equation}
\mathrm{d}\bar{\Omega}_{k}^{j}=\alpha_{l}^{j}\wedge\bar{\Omega}_{k}^{l}+\beta_{r}^{j}\wedge\Omega_{k}^{r}\label{eq:DerivedBasic}
\end{equation}
holds for suitable 1-forms $\alpha_{l}^{j},\beta_{r}^{j}$. 

If $\{S_{k}^{(1)},\mathrm{d}t\}$ is integrable (Frobenius theorem,
\cite{Bryant}), we have furthermore that
\begin{equation}
\mathrm{d}\bar{\Omega}_{k}^{j}=\gamma_{l}^{j}\wedge\bar{\Omega}_{k}^{l}+\rho^{j}\wedge\mathrm{d}t\label{eq:DerivedBasic-1}
\end{equation}
is met, for suitable 1-forms $\gamma_{l}^{j},\rho^{j}$. 
\begin{prop}
\label{prop:PropoDerived}Let us consider any sequence of Pfaffian
systems $\ldots\subset S_{2}\subset S_{1}\subset S_{0}$ with $S_{0}$
as in (\ref{eq:Pfaff1Org}) with elements $S_{k}$ together with their
first derived systems $S_{k}^{(1)}$. Then for all $v\in\mathcal{V}(S_{k})^{\perp}$
we have\end{prop}
\begin{enumerate}
\item $v\rfloor\mathrm{d}S_{k}^{(1)}\subset S_{k}^{(1)}$ is met if $\{S_{k}^{(1)},\mathrm{d}t\}$
is integrable.
\item $v\rfloor\mathrm{d}S_{k}^{(1)}\subset S_{k}$ 
\end{enumerate}
The proof of the first claim (a) follows by evaluating $v\rfloor\mathrm{d}\bar{\Omega}_{k}^{j}$
using (\ref{eq:DerivedBasic}) and (\ref{eq:DerivedBasic-1}) for
$v\in\mathcal{V}(S_{k})^{\perp}$ and the second (b) can be shown
by using (\ref{eq:DerivedBasic}) in a straightforward manner.

Systems that are exactly linearizable by static feedback meet $\{S_{k}^{(1)},\mathrm{d}t\}$
is integrable for every $k$, see \cite{Tilbury}. 
\begin{cor}
From Proposition \ref{prop:PropoDerived} it follows that for systems
$S_{0}$ that are exactly linearizable by static feedback the sequence
of the derived flags corresponds to the sequence as in Corollary \ref{cor:SequenceOrig}
and $\mathcal{F}_{k}$ corresponds to $\mathcal{V}(S_{k})^{\bot}$
which is integrable by construction. \end{cor}
\begin{rem}
The interesting case are of course examples that are not exactly linearizable
by static feedback, i.e. $\{S_{k}^{(1)},\mathrm{d}t\}$ are not integrable,
since then a different filtration has to be considered that may lead
to an implicit triangular form.
\end{rem}

\subsection{A constructive method to derive $\mathcal{F}_{k}$ and $S_{k+1}$\label{sub:A-constructive-method}}

If one is able to construct the sequence as in Corollary \ref{cor:SequenceOrig},
then one eventually ends up with the form (\ref{eq:TRiagPfaff}) (by
relabeling the coordinates) where in each step the involutive distribution
$\mathcal{F}_{l}$ has to be integrated, in order to derive the coordinate
transformation. Thus, in principle a constructive method that generates
the implicit triangular decomposition is stated. If then the rank
and dimension condition as in corollary \ref{cor:SequenceOrig} hold
the system is 0-flat/1-flat, but it should be stressed that this method
is only sufficient for flatness, and a failure does not in general
prove that a system is not flat.

However, the construction of $S_{k+1}\subset S_{k}$ such that an
involutive distribution $\mathcal{F}_{k}$ can be found that meets
$\mathcal{F}_{k}\subset\mathcal{V}(S_{k})^{\perp}$ as well as $\mathcal{F}_{k}\subset\mathcal{C}(S_{k+1})$
is a difficult task and leads in general to partial differential equations.
Furthermore, since the choice $\mathcal{F}_{k}\subset\mathcal{V}(S_{k})^{\perp}$
as well as of $S_{k+1}$ with $\mathcal{F}_{k}\subset\mathcal{C}(S_{k+1})$
is not unique in general (branching points may appear) it might be
necessary to iterate the construction of $\mathcal{F}_{k}$ and $S_{k+1}$
(see section \ref{sub:ExampleLevine}) - it should be noted that based
on a simple necessary condition candidates for $\mathcal{F}_{k}$
and $S_{k+1}$ are singled out as shown next.

We have to construct $\mathcal{F}_{k}\subset\mathcal{V}(S_{k})^{\bot}$
and $S_{k+1}\subset S_{k}$ such that $\mathcal{F}_{k}\rfloor\mathrm{d}S_{k+1}\subset S_{k+1}$
is met. Then also the necessary condition 
\begin{equation}
\mathcal{F}_{k}\rfloor\mathrm{d}S_{k+1}\subset S_{k}\label{eq:Ness1}
\end{equation}
holds, since $S_{k+1}\subset S_{k}$. For $\mathcal{V}(S_{k})^{\bot}=\{v_{i}\}$
and $S_{k}=\{\omega^{j}\}$ with $r=\mathrm{dim}(S_{k})$ we derive
the purely algebraic conditions (necessary conditions)
\begin{equation}
c^{i}v_{i}\rfloor\mathrm{d}(a_{j}\omega^{j})\wedge(\omega^{1}\wedge\ldots\wedge\omega^{r})=0\label{eq:Ness2}
\end{equation}
where $c^{i}$ and $a_{j}$ depend on all the system variables.

It should be noted that due to the requirement $\mathrm{dim}(S_{k})=\mathrm{dim}(S_{k+1})+\mathrm{dim}(\mathcal{F}_{k})$
one has to find $\mathrm{dim}(S_{k})-\mathrm{dim}(\mathcal{F}_{k})$
independent solutions $a_{j}\omega^{j}$ for $S_{k+1}$. From (b)
in Proposition \ref{prop:PropoDerived} we see that $S_{k}^{(1)}$
fulfills this necessary condition independently of $\mathcal{F}_{k}$.
If we furthermore assume that $S_{k}^{(1)}\subset S_{k+1}$ then the
construction of $S_{k+1}$ and $\mathcal{F}_{k}$ can be simplified
further as demonstrated in the next section in great detail. The strategy
is now to solve the necessary conditions (\ref{eq:Ness1}) or which
is the same (\ref{eq:Ness2}) and to generate solutions for which
then finally the criteria $\mathcal{F}_{k}\rfloor\mathrm{d}S_{k+1}\subset S_{k+1}$
has to be checked as well as the parametrization as in corollary \ref{cor:SequenceOrig}.

\section{Examples}

We now present two examples, in the first one we show how one algorithmically
can compute the sequence of codistributions using the necessary condition
(\ref{eq:Ness2}) and the second example demonstrates a case where
the algorithm stops in a dead end, and another iteration is at need.

\subsection{The motivating example revisited\label{sec:The-motivating-example} }

Let us write the equations (\ref{eq:ExSin}) as a Pfaffian system
of the form $S_{0}=\{\omega_{0}^{1},\omega_{0}^{2},\omega_{0}^{3}\}$
with 
\begin{eqnarray}
\omega_{0}^{1} & = & \mathrm{d}x{}^{1}-u^{1}\mathrm{d}t\nonumber \\
\omega_{0}^{2} & = & \mathrm{d}x{}^{2}-u^{2}\mathrm{d}t\label{eq:ExampleSinForm1}\\
\omega_{0}^{3} & = & \mathrm{d}x^{3}-\sin\left(\frac{u^{1}}{u^{2}}\right)\mathrm{d}t,\nonumber 
\end{eqnarray}
then we obtain the following proposition regarding the first reduction
step.
\begin{prop}
Given the system $S_{0}$ as in (\ref{eq:ExampleSinForm1}) we derive
a splitting of the form $S_{0}=S_{1}\oplus S_{1,c}$ as well as $v_{0}$
that meets $v_{0}\in\mathcal{V}(S_{0})^{\perp}$ and $v_{0}\in\mathcal{C}(S_{1})$.
Indeed, 
\[
v{}_{0}=u^{1}\partial_{u^{1}}+u^{2}\partial_{u^{2}}
\]
and $S_{1}=\{\omega_{1}^{1}=\omega_{0}^{3},\omega_{1}^{2}=u^{2}\omega_{0}^{1}-u^{1}\omega_{0}^{2}\}$
with 
\begin{equation}
\omega_{1}^{1}=\mathrm{d}x^{3}-\sin\left(\frac{u^{1}}{u^{2}}\right)\mathrm{d}t\,,\,\,\,\omega_{1}^{2}=u^{2}\mathrm{d}x^{1}-u^{1}\mathrm{d}x^{2}\label{eq:D1Example}
\end{equation}
as well as the complement $S_{1,c}=\{\omega_{1,c}^{3}\}$ with $\omega_{1,c}^{3}=\omega_{0}^{2}=\mathrm{d}x{}^{2}-u^{2}\mathrm{d}t$
possess the required properties.
\end{prop}
The proof of this proposition follows from the observation that $\mathcal{V}(S_{0})^{\bot}=\{\partial_{u^{1}},\partial_{u^{2}}\}$
and that $v_{0}\rfloor\mathrm{d}S_{1}\subset S_{1}$ as desired. We
will now show how one can derive $S_{1}$ and $v_{0}$. 

\begin{calc}The first derived system $S_{0}^{(1)}$ is given by the
single form
\[
\Phi=\cos(\frac{u^{1}}{u^{2}})\left(\frac{u^{1}}{u^{2}}\mathrm{d}x^{2}-\mathrm{d}x^{1}\right)+\mathrm{u^{2}(d}x^{3}-\sin(\frac{u^{1}}{u^{2}})\mathrm{d}t)
\]
and a basis for $S_{0}$ can be alternatively given as $S_{0}=\{\omega_{0}^{1},\omega_{0}^{2},\Phi\}$.
To construct $S_{1}$ we assume that $S_{0}^{(1)}\subset S_{1}$ and
consider the relation (according to (\ref{eq:Ness2})) 
\begin{equation}
(c_{1}^{1}\partial_{u^{1}}+c_{1}^{2}\partial_{u^{2}})\rfloor\mathrm{d}(a_{1}^{1}\mathrm{\omega_{0}^{1}+a_{2}^{1}}\omega_{0}^{2}+a_{3}^{1}\Phi)\wedge\omega_{0}^{1}\wedge\omega_{0}^{2}\wedge\Phi=0\label{eq:D1formula}
\end{equation}
where $c_{1}^{i}$ and $a_{i}^{1}$ are functions of all system variables
that have to be computed. From (\ref{eq:D1formula}) we are left with
the equation $c_{1}^{1}a_{1}^{1}+c_{1}^{2}a_{2}^{1}=0$ or $a_{1}^{1}=-a_{2}^{1}\frac{c_{1}^{2}}{c_{1}^{1}}$.
This means that the forms $\Phi$ and $\omega_{0}^{2}-\frac{c_{1}^{2}}{c_{1}^{1}}\omega_{0}^{1}$
fulfill the necessary conditions for the vector field $v_{0}=c_{1}^{1}\partial_{u^{1}}+c_{1}^{2}\partial_{u^{2}}$.
To determine $c_{1}^{1}$ and $c_{1}^{2}$ we consider the criteria
$v_{0}\rfloor\mathrm{d}S_{1}\subset S_{1}$ and we derive the relation
\begin{equation}
\left((c_{1}^{1}\partial_{u^{1}}+c_{1}^{2}\partial_{u^{2}})\rfloor(\mathrm{d}\Phi)\right)\wedge(\omega_{0}^{2}-\frac{c_{1}^{2}}{c_{1}^{1}}\omega_{0}^{1})\wedge\Phi=0.\label{eq:DFrel1}
\end{equation}
For the solution of (\ref{eq:DFrel1}) of the form $c_{1}^{1}=c_{1}^{2}\frac{u^{1}}{u^{2}}$
we have that $\Phi$ and $\omega_{0}^{2}-\frac{u^{2}}{u^{1}}\omega_{0}^{1}$
clearly correspond to $S_{1}$ as in (\ref{eq:D1Example}) as can
be checked easily (by linear combinations) and that 
\begin{equation}
(v_{0}\rfloor\mathrm{d}\omega_{1}^{2})\wedge\omega_{1}^{1}\wedge\omega_{1}^{2}=0\label{eq:DFrel2}
\end{equation}
is fulfilled, such that $v_{0}\rfloor\mathrm{d}S_{1}\subset S_{1}$
is met, because of (\ref{eq:DFrel1},\ref{eq:DFrel2}). 

\end{calc}To straighten out $v_{0}$ we consider the coordinate transformation
$(x^{1},x^{2},x^{3},u^{1},u^{2})=\varphi_{0}(w^{1},w^{2},w^{3},w^{4},\hat{w})$
with $x^{i}=w^{i}$ for $i=1,2,3$ and 
\begin{eqnarray*}
u^{1} & = & e^{\hat{w}}w^{4}\\
u^{2} & = & e^{\hat{w}}
\end{eqnarray*}
which is based on the flow of $v_{0}.$ In new coordinates we obtain
a basis for $S_{1}$ as
\begin{equation}
\begin{array}{cclcc}
\omega_{1}^{1} & = & \mathrm{\mathrm{d}}w^{3}-\sin\left(w^{4}\right)\mathrm{\mathrm{d}}t\\
\omega_{1}^{2} & = & \mathrm{d}w^{1}-w^{4}\mathrm{d}w{}^{2}
\end{array}\label{eq:S1ex}
\end{equation}
and for the complement $S_{1,c}=\{\omega_{1,c}^{3}\}$ with $\omega_{1,c}^{3}=\mathrm{d}w^{2}-e^{\hat{w}}\mathrm{d}t$
and it can be checked easily that $\mathrm{dim}(\mathcal{F}_{0})=1$,
$\mathrm{dim}(S_{1})+1=\mathrm{dim}(S_{0})$ and that the Jacobian
$\partial_{\hat{w}}(\dot{w}_{2}-e^{\hat{w}})$ has maximal rank. 
\begin{rem}
We want to point out again, that $v_{0}$ is a Cauchy characteristic
vector field for $S_{1}$, i.e. $v_{0}\in\mathcal{C}(S_{1})$ and
this guarantees that there is a basis for $S_{1}$ which does not
depend on the coordinate $\hat{w}$, since in new coordinates $\partial_{\hat{w}}\in\mathcal{C}(S_{1})$
is met.
\end{rem}
Then we continue our considerations with $S_{1}$ and the following
proposition states the second reduction step.
\begin{prop}
Given the system $S_{1}$ as in (\ref{eq:S1ex}) we derive a splitting
of the form $S_{1}=S_{2}\oplus S_{2,c}$ as well as $v_{1}$ that
meets $v_{1}\in\mathcal{V}(S_{1})^{\perp}$ and $v_{1}\in\mathcal{C}(S_{2})$
with 
\begin{eqnarray}
v_{1} & = & w^{4}\partial_{w^{1}}+\partial_{w^{2}}\label{eq:Vec2}
\end{eqnarray}
and 
\begin{equation}
S_{2}=\{\omega_{2}^{1}\}\,,\,\,\,\,\omega_{2}^{1}=\mathrm{d}w^{3}-\sin\left(w^{4}\right)\mathrm{d}t\label{eq:D2ex}
\end{equation}
and $S_{2,c}=\{\omega_{2,c}^{2}\}$ with $\omega_{2,c}^{2}=\mathrm{d}w^{1}-w^{4}\mathrm{d}w{}^{2}$.
\end{prop}
The proof follows again from the fact that $\mathcal{V}(S_{1})^{\perp}=\{w^{4}\partial_{w^{1}}+\partial_{w^{2}},\partial_{w^{4}}\}$
and $v_{1}\rfloor\mathrm{d}S_{2}\subset S_{2}.$ The construction
of $S_{2}$ can be performed in the same manner as above. (Observe
however that $S_{1}^{(1)}$ is empty, but from
\[
(c_{2}^{1}(w^{4}\partial_{w^{1}}+\partial_{w^{2}})+c_{2}^{2}\partial_{w^{4}})\rfloor\mathrm{d}(a_{1}^{2}\mathrm{\omega_{1}^{1}+a_{2}^{2}}\omega_{1}^{2})\wedge\omega_{1}^{1}\wedge\omega_{1}^{2}=0
\]
 that result follows at once). 

Based on the flow of $v_{1}$ we derive the map $w=\varphi_{1}(q,\hat{q})$
in the form
\[
\begin{array}{cllcccc}
w^{1} & = & \hat{q}q^{4} & , & w^{3} & = & q^{3}\\
w^{2} & = & \hat{q}+q^{2} & , & w^{4} & = & q^{4}.
\end{array}
\]
With 
\begin{equation}
(y^{1}=q^{3},\, y^{2}=q^{2},\,\hat{z}^{2}=q^{4},\,\hat{z}^{3}=\hat{q},\,\hat{z}^{4}=\hat{w})\label{eq:ChangeName}
\end{equation}
it is easily seen that the composition of $(x,u)=\varphi_{0}(w,\hat{w})$
and $w=\varphi_{1}(q,\hat{q})$ together with (\ref{eq:ChangeName})
gives the desired transformation $(x,u)=\varphi(z)$ as in (\ref{eq:Trafo1}).
Furthermore, the sequence of systems $S_{2}\subset S_{1}\subset S_{0}$
leads at once to the desired normal-form 
\begin{eqnarray*}
\omega_{d,0}^{1} & = & \mathrm{d}y^{1}-\sin\left(\hat{z}^{2}\right)\mathrm{d}t\\
\omega_{d,0}^{2} & = & \hat{z}^{3}\mathrm{d}\hat{z}^{2}-\hat{z}^{2}\mathrm{d}y^{2}\\
\omega_{d,0}^{3} & = & \mathrm{d}\hat{z}^{3}+\mathrm{d}y^{2}-e^{\hat{z}^{4}}\mathrm{d}t.
\end{eqnarray*}
as in (\ref{eq:SinNormal}). The flow parameters $\hat{w}$ and $\hat{q}$
correspond to the non-derivative variables $\hat{z}^{4}$ and $\hat{z}^{3}$,
respectively. Furthermore, $y^{2}$ is a flat output since $\mbox{\ensuremath{\omega_{d,0}^{2}}}$
is parameterizable with respect to $\hat{z}^{3}$ and $\partial_{q^{2}}\subset\mathcal{C}(S_{2})$
with $q^{2}=y^{2}$.

\subsection{A further example\label{sub:ExampleLevine}}

Let us consider the system $S_{0}=\{\omega_{0}^{1},\omega_{0}^{2},\omega_{0}^{3},\omega_{0}^{4}\}$
also treated in \cite{CharletSIAM} in a different context 
\begin{eqnarray*}
\omega_{0}^{1} & = & \mathrm{d}x^{1}-(x^{2}+x^{3}u^{2})\mathrm{d}t\\
\omega_{0}^{2} & = & \mathrm{d}x^{2}-(x^{3}+x^{1}u^{2})\mathrm{d}t\\
\omega_{0}^{3} & = & \mathrm{d}x^{3}-(u^{1}+x^{2}u^{2})\mathrm{d}t\\
\omega_{0}^{4} & = & \mathrm{d}x^{4}-u^{2}\mathrm{d}t
\end{eqnarray*}
where we again have $\mathcal{V}(S_{0})^{\bot}=\{\partial_{u^{1}},\partial_{u^{2}}\}$.
The triangular form is based on the decompositions $S_{0}=S_{1}\oplus S_{1,c}$
with $S_{1}=\{\omega_{1}^{1}=\omega_{0}^{1},\omega_{1}^{2}=\omega_{0}^{2},\omega_{1}^{3}=\omega_{0}^{4}\}$
\begin{eqnarray*}
\omega_{1}^{1} & = & \mathrm{d}x^{1}-(x^{2}+x^{3}u^{2})\mathrm{d}t\\
\omega_{1}^{2} & = & \mathrm{d}x^{2}-(x^{3}+x^{1}u^{2})\mathrm{d}t\\
\omega_{1}^{3} & = & \mathrm{d}x^{4}-u^{2}\mathrm{d}t
\end{eqnarray*}
and $S_{1,c}=\{\omega_{1,c}^{4}=\mathrm{d}x^{3}-(u^{1}+x^{2}u^{2})\mathrm{d}t\}$
as well as on $S_{1}=S_{2}\oplus S_{2,c}$ with $S_{2}=\{\omega_{2}^{1}=\omega_{1}^{1}-u^{2}\omega_{1}^{2},\,\omega_{2}^{2}=\omega_{1}^{3}\}$
with 
\begin{eqnarray*}
\omega_{2}^{1} & = & \mathrm{d}x^{1}-u^{2}\mathrm{d}x^{2}-(x^{2}-x^{1}(u^{2})^{2})\mathrm{d}t\\
\omega_{2}^{2} & = & \mathrm{d}x^{4}-u^{2}\mathrm{d}t
\end{eqnarray*}
and $S_{2,c}=\{\omega_{2,c}^{3}=\mathrm{d}x^{2}-(x^{3}+x^{1}u^{2})\mathrm{d}t\}$.
The distributions $\mathcal{F}_{0}=\{\partial_{u^{1}}\}\subset\mathcal{V}(S_{0})^{\perp},\mathcal{\, F}_{1}=\{\partial_{x^{3}}\}\subset\mathcal{V}(S_{1})^{\perp}$
and $\mathcal{F}_{2}=\{\partial_{x^{2}}+u^{2}\partial_{x^{1}}\}\subset\mathcal{V}(S_{2})^{\perp}$
were used and the flat outputs $y^{1}=x^{1}-u^{2}x^{2}$ and $y^{2}=x^{4}$
follow at once by applying a coordinate transformation based on the
flow of $\mathcal{F}_{2}$ and regarding $\mathcal{F}_{0}$ and $\mathcal{F}_{1}$
no coordinate transformation is at need, since $u^{1}$ and $x^{3}$
are already non-derivative variables.
\begin{rem}
Also in this example we have that $S_{0}^{(1)}\subset S_{1}$ and
$S_{1}^{(1)}\subset S_{2}$ which enables one to construct the solutions
based on the necessary condition (\ref{eq:Ness2}) very easily. 
\end{rem}
However, a different possible solution for $S_{1}\oplus S_{1,c}$
(branching point) can be based on choosing the distribution $\mathcal{F}_{0}=\{\partial_{u^{1}},\partial_{u^{2}}\}$
together with $S_{1}=S_{0}^{(1)}=\{\omega_{1}^{1},\omega_{1}^{2}\}$
\begin{eqnarray*}
\omega_{1}^{1} & = & \mathrm{d}x^{1}-x^{3}\mathrm{d}x^{4}-x^{2}\mathrm{d}t\\
\omega_{1}^{2} & = & \mathrm{d}x^{2}-x^{1}\mathrm{d}x^{4}-x^{3}\mathrm{d}t
\end{eqnarray*}
and $S_{1,c}=\{\mathrm{d}x^{3}-(u^{1}+x^{2}u^{2})\mathrm{d}t\,,\mathrm{d}x^{4}-u^{2}\mathrm{d}t\}$,
where obviously $S_{1,c}$ is parameterizable with respect to $u^{1}$
and $u^{2}$. This choice for $\mathcal{F}_{0}$ and $S_{1}$ however
leads to a 'dead end' since for $S_{1}$ the necessary condition (\ref{eq:Ness2})
does not lead to a splitting $S_{1}=S_{2}\oplus S_{2,c}$.

\section{Discussion}

We have characterized a suitable normal form for 1-flat systems, which
is in implicit triangular shape, see (\ref{eq:TRiagPfaff}), that
possesses the properties as in proposition \ref{prop:Properties}
based on exterior algebra. Furthermore, we have discussed a constructive
calculation scheme to transform 1-flat systems into that desired form.
It should be mentioned again that we only provide sufficient conditions
for a system to be 1-flat and that the constructive algorithm is in
general not unique, and iterations might be necessary. Nevertheless,
we believe that the presented normal-form is of interest in the analysis
of the flatness problem, and our examples show that this implicit
triangular form can be achieved by successive coordinate transformations
in a rather straightforward manner. Additionally, the well known Brunovsky
form for systems that are linearizable by static feedback is naturally
included in our approach, based on proposition \ref{prop:PropoDerived}. 

\section*{Acknowledgment}

Markus Sch\"{o}berl is an APART fellowship holder of the Austrian Academy of Sciences. 
\bibliographystyle{plain}   
\bibliography{./maxibib2}

                                        % in the appendices.
\end{document}